\documentclass[]{interact}
\usepackage[usenames,dvipsnames]{xcolor}
\usepackage{latexsym,bezier,enumerate,dcolumn}
\usepackage{amsmath,amsthm,amsopn,amstext,amscd,amsfonts,amssymb}
\usepackage{graphics,lscape,epsfig,fancyhdr}
\usepackage{graphicx,comment,url,times,subfigure}
\usepackage{relsize}
\usepackage{cancel}
\usepackage{marginnote}
\usepackage[highstructure]{accessibility}
\usepackage[title]{appendix}
\usepackage{setspace}
\usepackage{color}
\usepackage{longtable}
\usepackage{algorithm}

\usepackage{natbib}
\bibpunct[, ]{(}{)}{;}{a}{}{,}
\renewcommand\bibfont{\fontsize{10}{12}\selectfont}

\theoremstyle{plain}

\theoremstyle{definition}

\newcommand{\BB}{{\cal B}}
\newcommand{\C}{{\cal C}}

\newcommand{\I}{{\cal I}}
\newcommand{\J}{{\cal J}}
\newcommand{\K}{{\cal K}}

\newcommand{\M}{{\cal M}}
\newcommand{\N}{{\cal N}}

\renewcommand{\S}{{\cal S}}

\newcommand{\bea}{\begin{eqnarray}}
\newcommand{\eea}{\end{eqnarray}}

\begin{document}


\title{A methodology for the  Cross-dock Door Platforms design under
  uncertainty}
\author{
\name{Laureano
  F. Escudero\textsuperscript{a}\thanks{\textsuperscript{a} Email: laureano.escudero@urjc.es}
 and
  M. Araceli Gar\'in\textsuperscript{b}\thanks{\textsuperscript{b} Corresponding
    author email: mariaaraceli.garin@ehu.eus}
 and
 Aitziber Unzueta\textsuperscript{c}\thanks{\textsuperscript{c}Email:  aitziber.unzueta@ehu.eus}
}
\affil{\textsuperscript{a} Statistics and
    Operations Research Area. Universidad Rey Juan Carlos, URJC,
    M\'ostoles (Madrid), Spain\break
     \textsuperscript{b} Quantitative Methods Department.
                 Universidad del Pa\'{\i}s Vasco, UPV/EHU, Bilbao
                 (Bizkaia), Spain\break
     \textsuperscript{c} Applied Mathematics Department.
                Universidad del Pa\'{\i}s Vasco, UPV/EHU, Bilbao (Bizkaia), Spain}
}

\maketitle


\noindent\textbf{Abstract.}
The Cross-dock Door Design Problem (CDDP) consists of deciding on the
number and capacity  of inbound  and outbound doors for receiving
product pallets from origin nodes and exiting them to destination
nodes. The uncertainty, realized in scenarios, lies in the occurrence
of these nodes, the number and cost of the pallets, and the capacity’s
disruption of the doors. The CDDP is represented using a stochastic
two-stage binary quadratic model (BQM). The first stage decisions are
related to the cross-dock infrastructure design, and the second stage
decisions are related to the node-to-door’s assignments. This is the
first time, as far as we know, that a stochastic two-stage BQM has
been presented for minimizing the construction cost of the
infrastructure and its exploitation expected cost in the
scenarios. Given the difficulty of solving this combinatorial problem,
a mathematically equivalent MILP formulation is introduced. However,
searching an optimal solution is  still impractical for commercial
solvers. Thus, a  scenario cluster decomposition-based  matheuristic
algorithm is introduced to obtain feasible solutions with  small
optimality gap and reasonable computational effort.
A broad study to validate the proposal gives solutions with a much
smaller gap than the ones provided by a state-of-the-art general
solver. In fact,  the proposal provides solutions with a 1 to 5\%
optimality gap, while the solver does it  with up to a 12\% gap, if
any, and requires a wall time two orders of magnitude higher.


\section{Introduction and motivation}\label{sec:intro}
\begin{figure}[h]
\begin{center}
   \includegraphics[height=4.cm]{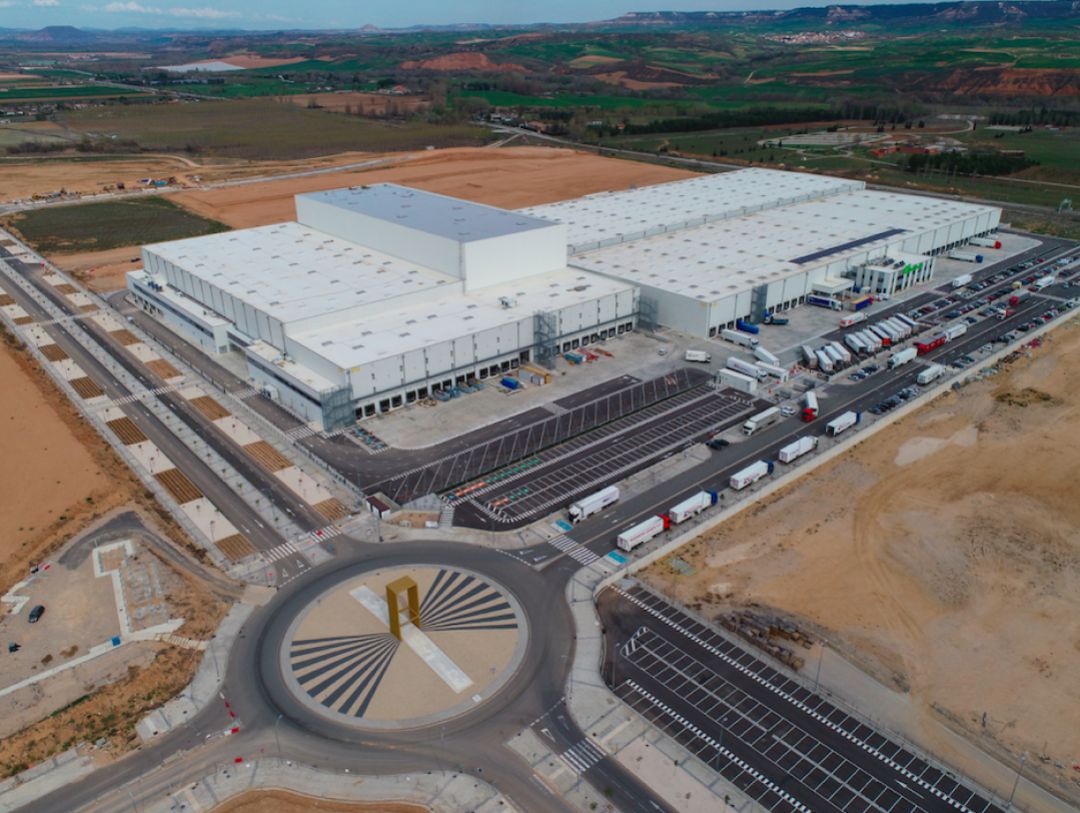}
\caption{Cross-docking door terminal in Guadalajara (Spain). Source: Alimarket}\label{fig:cross-dock-Guada}
\end{center}
\end{figure}

\begin{figure}[h]
\begin{center}
\includegraphics[width=8.cm]{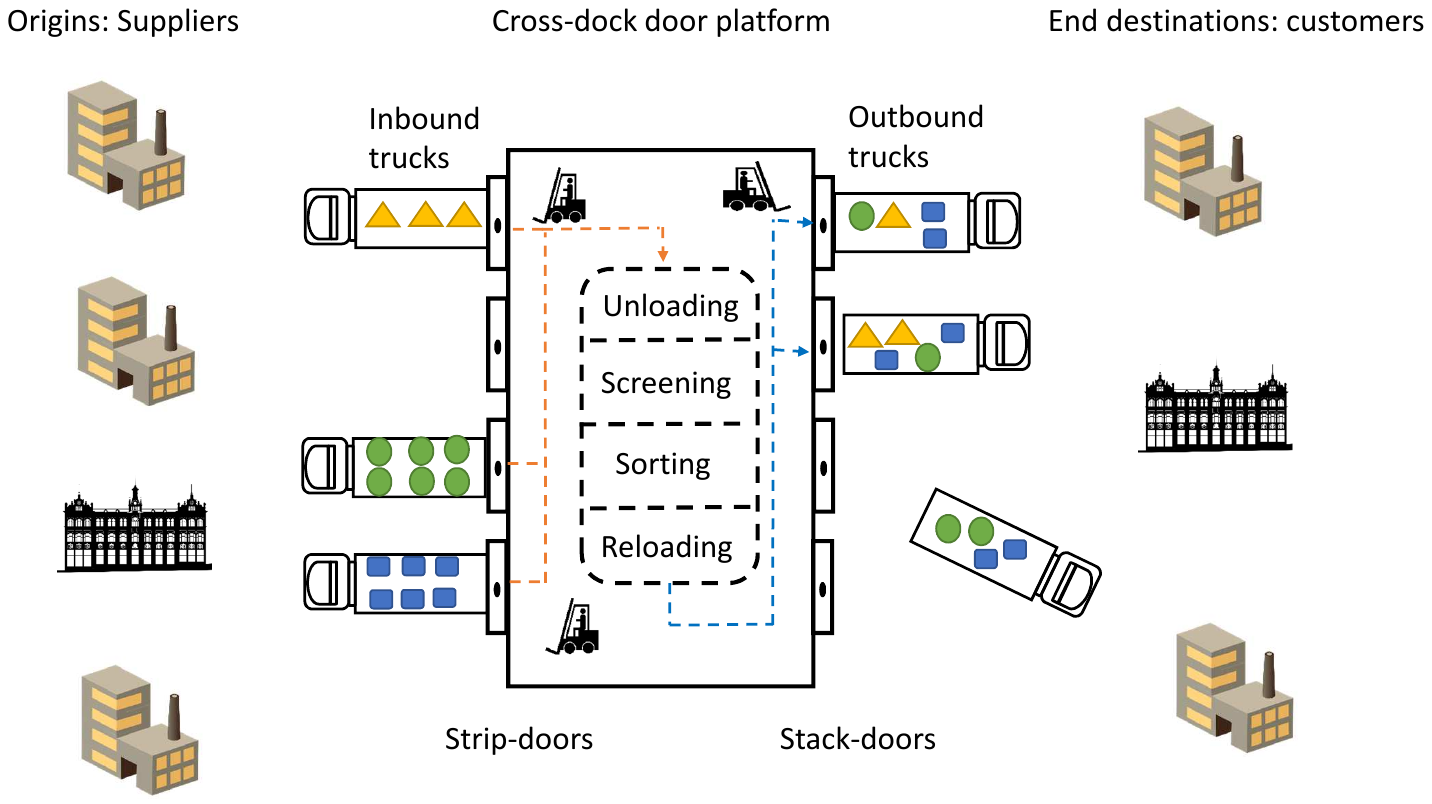}
\vskip -1.5cm
\caption{Cross-dock infrastructure system}\label{fig:cross-dock}
\end{center}
\end{figure}

Given a network with a set of supplying (i.e. origin) nodes for different product types and a set of receiving (i.e. destination) nodes for these products, usually in smaller quantities, a cross-dock entity may serve as a consolidation point.
The origin nodes can deliver the material at the cross-dock so that, after being classified by type and destination, it can be transported to the destination nodes.
A cross-dock infrastructure has a number of inbound doors, which are called
strip doors, and a number of outbound doors, which are called stack doors, and each of them
has a capacity for pallet handling during a given time period.
Fig. \ref{fig:cross-dock-Guada} depicts an image of a real terminal in Guadalajara (Spain).
The classical operation process of a cross-dock is depicted in
Fig. \ref{fig:cross-dock}.

Two main types of optimization problems arise in cross-dock dealing. One of them is the Cross-dock Door Design Problem (CDDP), the subject of this work,
where the number of strip and stack doors is not given, nor the related nominal capacities; additionally,
the set of origin and destinations nodes is uncertain, along with the product volume to be served.
The other main type of optimization problem is the deterministic Cross-dock Door Assignment Problem (CDAP),
where the number of strip and stack doors is given as well as the related capacity.
The origin and destination nodes are also known, along with the product volume to be served.
The CDAP is also a difficult binary quadratic problem (BQP) whose constraint system is composed of two Generalised Assignment Problems (GAPs);
one is related to the strip doors and the other is related to the stack doors to which the origin nodes and destination nodes are to be assigned, respectively.
These two types of problems are linked by the quadratic objective function.
The scientific treatment of cross-docking (CD) problems is a young discipline;
in fact, most of the literature has been published during the last 25 years.

Different approaches for CDAP solving and comprehensive reviews can be found in \cite{VanBelle12} and \cite{Goodarzi20}, mainly for metaheuristic algorithms.
In particular, see \cite{Guignard12,kuzu17,Nassief19,Gelareh20,Goodarzi20,Guignard20,cdap23}; and \cite{Li24}.
Comprehensive reviews on uncertainty have also been recently published by \cite{ArdakaniFei20} and \cite{BuakumWisittipanich19}.
In particular, see \cite{Goodarzi20} for uncertainty related to disruption, reliability, and reallocation issues and \cite{Xi20} for uncertainty related to arrival and operational times; the former studies the cross-dock truck scheduling problem, which provides an exact algorithm for column and constraint generation.
Moreover, the literature on the CDDP is very scarce.
The cross-dock design is usually considered within the design of a supply network, which is planned along a time horizon under uncertainty,
where the details of each individual door assignment are not taken into account.
Some works that are representative of the state of the art for this problem are \cite{Mousavi14,Soanpet12,Gallo22} and \cite{TaheriTaft24}.
In particular, the first two works consider the locations of multiple cross-dock centres, the assignment of origin and destination nodes to the centres, receiving, stocking, and exiting products, and vehicle routing scheduling for supplying and delivering products along a time horizon.
An MILP fuzzy optimization method is proposed in the first work, and an MILP robust optimization method is proposed in the second work to cope with uncertain critical parameters, such as product supplies and demands, the capacity of vehicles, the transportation time needed for each vehicle to move between nodes, and transportation and  operational costs.
\cite{Soanpet12} presents two stochastic models for CDDP in a multi-product supply network planning problem for given sets of origin and destination nodes. The cross-dock centres are chosen to consolidate the products and lower transportation and handling costs.
However, the individual doors, the capacity, and product assignment are not considered.
The first model is a chance constraint stochastic integer quadratic model used to locate a given number of centres,
where each one has a pre-defined capacity, as well as to assign nodes and vehicles to the centres.
The uncertainty lies in the capacity of each centre; it is assumed to follow the normal distribution.
The goal is to minimize the expected total cost, provided that the locations of the centres, as well as the transportation and handling operations, allow the probability of satisfying the capacity of each centre to be larger than a given threshold.
The second model is a two-stage stochastic integer model,
where the first stage deals with the infrastructure of the centres.
The second stage deals with the supply network operations (i.e. assigning nodes and vehicles to the centres) for each scenario,
where the uncertainty lies in the varying capacity disruption of the centres.
The goal is to minimize the expected total cost of the centres' locations and product routing and handling.
\cite{Gallo22} present a genetic metaheuristic algorithm for solving a two-stage stochastic MILP model in the cross-dock infrastructure. The goal is to minimize the earliness and tardiness penalty costs of delivering products from inbound doors when the intended delivery time is within the time windows for exporting these products through outbound doors under uncertain truck arrivals.
A computational experiment is performed using a real-world pilot case.
\cite{TaheriTaft24} present a multi-objective robust MILP approach for minimizing the loading/unloading transportation costs, as well as the product distribution and shipping durations, while scheduling and routing vehicles through a multi cross-dock-system,
taking into account the uncertainties in demand.

\subsection{CDDP description}
The pure CDDP consists of deciding on the number of strip and stack doors and their nominal capacity for receiving product pallets through the strip doors from the origin nodes (i.e. the suppliers), consolidating the products in a collection of mixed-destination pallets, and finally, delivering them through the stack doors to the destination nodes.
The uncertainty lies in the set of origin and destination nodes and in the number of product pallets that must traverse the cross-dock from the origin nodes to the destination ones in a given time period.
The other source of uncertainty that is considered is the capacity's disruption of the doors, due to sabotage, misuse, etc.
The uncertainty is realized in a finite set of scenarios, so that it represents the occurrence of the CDDP elements within the time horizon.
The goal is to minimize the cross-dock infrastructure construction cost plus the expected cost of the assignments of the nodes to the doors in the scenarios, subject to the related net capacity constraints.
The literature is very scarce for this type of cross-dock problem; see \cite{Shams-Shemirani23}, where heuristic simulation-related algorithms are presented to evaluate different historical scenarios to choose the number of doors and loaders in a cross-dock infrastructure design.

\subsection{Contributions}
\begin{description}
\item (i) Introduce a two-stage stochastic BQP model for cross-dock door design planning (CDDP-TS)
    as an extension of the deterministic CDAP to deal with the uncertainty.
    The first stage is devoted to the strategic decisions (i.e. the number of strip and stack doors and related nominal capacities),
    before the uncertainties related to the set of origin and destination nodes and the disruption of the doors capacity are unveiled.
    The second stage is devoted to the operational decisions (i.e. the assignment of doors to origin and destination nodes in the scenarios).
    The objective function to minimize is composed of the binary linear cross-dock door infrastructure construction cost and the binary quadratic cost function related to the CDAP scenario node-to-door assignment.
    See Section \ref{sec:CDDP-TS}.
\item (ii) Develop a Linearized mixed Integer Programming problem (LIP) as a reformulation that is mathematically equivalent to the quadratic one
    by using the Reformulation Linearization Technique (RLT1), see \cite{SheraliAdams94};
    as a result, there is a high number of additional continuous variables.
    See Section \ref{sec:LIP}.
\item (iii) Introduce two options for a scenario Cluster Decomposition (CD) of model LIP,
    where the special structure of the problem is benefitted from.
    Option 1 decomposes model LIP into scenario clusters for the first stage variables.
    Option 2 additionally does this for the strip and stack door types in each cluster.
    See Section \ref{sec:cd}.
\item (iv) Develop a matheuristic algorithm, based on feasible first stage solutions of the CD model,
    to obtain feasible solutions for the original model CDDP-TS.
    This algorithm considers a linear search approach for large-scale instances to exploit the scenario structure of the second stage submodels.
    See Section \ref{sec:SCS4B}.
\end{description}

The main results of a broad computational experiment are reported in Section \ref{sec:results} to validate the proposal made in this work.
The proposed method provides solutions with a 1 to 5\% optimality gap, while a state-of-the-art solver provides solutions with up to a 12\% gap and requires a wall time that is two orders of magnitude higher than the wall time required by the proposed method, in the case in which it provides a solution.
Section \ref{sec:con} draws some conclusions.

\section{A two-stage binary quadratic model for Cross-dock Door Design under uncertainty}\label{sec:CDDP-TS}
Let us introduce the binary quadratic model for the CDDP under uncertainty, as an extension of the deterministic CDAP; see \cite{Guignard12,Gelareh20,Guignard20}; and \cite{cdap23}, among others.

\noindent\textit{First stage sets for the cross-dock infrastructure:}
\begin{description}
\item[$\I$,] candidate strip doors, without including the inbound outsourcing `door' $i=0$.
\item[$\K_i$,] capacity level for candidate strip door $i$, for $i\in\I$,
    such that $\K_i=\{0,1,...,|\K_i|\}$.
\item[$\J$,] candidate stack doors, without including the outbound outsourcing `door' $j=0$.
\item[$\K_j$,] capacity level for candidate stack door $j$, for $j\in\J$,
    such that $\K_j=\{0,1,...,|\K_j|\}$.
\end{description}

\noindent\textit{First stage parameters:}
\begin{description}
\item[$E_{ij}$,] distance between strip door $i$ and stack door $j$, for $i\in\I, \, j\in\J$.
\item[$S_{ki}$,] nominal capacity of alternative level $k$ in strip door $i$,
    where $S_{ki} < S_{k'i}$ for $k,k'=1,2,...,|\K_i|: \, k<k'$, $i\in\I$.

Let $S_{0i}$ denote the so-called \textit{basic capacity} of strip door $i$.
It is a special capacity that helps to avoid outsourcing solutions (see below).

\item[$R_{kj}$,] nominal capacity of alternative level $k$ in stack door $j$,
    where $R_{kj} < R_{k'j}$ for $k,k'=1,2,...,|\K_j|: \, k<k'$, $j\in\J$.

Let $R_{0j}$ denote the \textit{basic capacity} of stack door $j$.
It is a special capacity that helps to avoid outsourcing solutions (see below).

\item[$\overline{I}$] and $\overline{J}$, upper bounds on the number of strip and stack doors, respectively, that are allowed in the cross-dock.
\end{description}

\noindent\textit{Strategic (first stage) costs of cross-dock door building:}
\begin{description}
\item[$F_{ki}$,] cost of installing capacity level $k$ in strip door $i$, for $k\in\K_i, \, i\in\I$,
    with a high enough cost for $k=0$.

\item[$F_{kj}$,] cost of installing capacity level $k$ in stack door $j$, for $k\in\K_j, \, j\in\J$,
    with a high enough cost for $k=0$.

\item[$F_0$,] a high enough penalization for considering outsourcing `doors' $i=0$ and $j=0$.
\end{description}

\noindent\textit{Second stage sets for origin and destination nodes to be served in the cross-dock under uncertainty:}
\begin{description}
\item[$\Omega$,] scenarios.
\item[$\M^\omega$,] origin nodes in scenario $\omega$.
\item[$\N^\omega$,] destination nodes in scenario $\omega$.
\end{description}

\noindent\textit{Second stage parameters for origin node $m$ and destination node $n$ under scenario $\omega$,
for $m\in\M^\omega, \, n\in\N^\omega, \, \omega\in\Omega$:}
\begin{description}
\item[$w^\omega$,] weight factor representing the likelihood that is associated with the scenario,
    so that $\sum_{\omega\in\Omega}w^\omega=1$.
\item[$H_{mn}^\omega$,] commodity volume (i.e. number of pallets) to be consolidated in the cross-dock from node $m$ to node $n$.
\item[$S_m^\omega$,] total commodity volume (i.e. number of product pallets) to enter the cross-dock through any strip door from node $m$,
    so that $S_m^\omega = \sum_{n\in\N^\omega}H_{mn}^\omega$.
\item[$R_n^\omega$,] total commodity volume (i.e. number of pallets) to exit the cross-dock through any stack door to node $n$,
    so that $R_n^\omega = \sum_{m\in\M^\omega}H_{mn}^\omega$.
\item[$D_i^\omega$,] fraction of disruption of the capacity of strip door $i$, for $i\in\I$.
\item[$D_j^\omega$,] fraction of disruption of the capacity of stack door $j$, for $j\in\J$.
\end{description}

\noindent\textit{Updated door set for origin-destination nodes under scenario $\omega$, for $\omega\in\Omega$:}
\begin{description}\parskip 0.5mm
\item[$\I_m^\omega\subseteq\I$,] the subset $\{i\}$ of strip doors expressed as\\
      $\{ i\in\I: \exists k\in\K_i \text{, where } S_m^\omega \leq (1-D_i^\omega)S_{ki} \}$, for $m\in\M^\omega$.
\item[$\J_n^\omega\subseteq\J$,] the subset $\{j\}$ of stack doors expressed as\\
      $ \{ j\in\J: \exists k\in\K_j \text{, where } R_n^\omega \leq (1-D_j^\omega)R_{kj} \}$, for $n\in\N^\omega$.
\end{description}

\noindent\textit{Operational (second stage) costs of handling the commodity volume $H_{mn}^\omega$ from origin node $m$ in the cross-dock and exiting it to destination node $n$ under scenario $\omega$, for $m\in\M^\omega, \, n\in\N^\omega, \, \omega\in\Omega$:}
\begin{description}
\item[$G_{minj}^\omega$,] cost of the standard operation on the commodity's pallets entering the cross-dock from origin node $m$ in the cross-dock through strip door $i$ and exiting it through stack door $j$ to destination node $n$, based on $E_{ij} H_{mn}^\omega$, for $i\in\I, \, j\in\J$.

\item[$G_{m0n0}^\omega$,] a high enough penalization for considering the outsourcing `doors', for $i=0$ and $j=0$.
\end{description}

\noindent\textit{First-stage binary variables for cross-dock infrastructure design planning:}
\begin{description}
\item[$\alpha_{ki} = 1$,] if capacity level $k$ is installed in strip door $i$; otherwise, it is 0, for $k\in\K_i, \, i\in\I$.
\item[$\beta_{kj}  = 1$,] if capacity level $k$ is installed in stack door $j$; otherwise, it is 0, for $k\in\K_j, \, j\in\J$.
\end{description}

\noindent\textit{Second-stage binary variables under scenario $\omega$,
for $m\in\M^\omega, \, n\in\N^\omega, \, \omega\in\Omega$:}
\begin{description}
\item[$\alpha_0^\omega = 1$,] if inbound outsourcing `door' $i=0$ is considered; otherwise, it is 0.
\item[$\beta_0^\omega  = 1$,] if outbound outsourcing `door' $j=0$ is considered; otherwise, it is 0.
\item[$x_{mi}^\omega   = 1$,] if origin node $m$ is assigned to strip door $i$; otherwise, it is 0, for $i\in\I$.
\item[$x_{m0}^\omega   = 1$,] if origin node $m$ is \textit{not} assigned to any strip door; otherwise, it is 0.
\item[$y_{nj}^\omega   = 1$,] if destination node $n$ is assigned to stack door $j$; otherwise, it is 0, for $j\in\J$.
\item[$y_{m0}^\omega   = 1$,] if destination node $n$ is \textit{not} assigned to any stack door; otherwise, it is 0.
\end{description}
Note that if the solution retrieved from model CDDP-TS is such that $\alpha_{0}^\omega = 1$ or $\beta_{0}^\omega = 1$,
then it is called an outsourcing solution (i.e. let us say an infeasible one).

The \noindent\textit{binary quadratic model CDDP-TS} can be expressed
\begin{subequations}\label{CDDP-obj}
\bea
& \label{Qa} z^* \, = \, \min \quad \displaystyle
       \sum_{i\in\I}\sum_{k\in\K_i}F_{ki} \alpha_{ki} + \sum_{j\in\J}\sum_{k\in\K_j}F_{kj} \beta_{kj}
     + \sum_{\omega\in\Omega}w^\omega F_0 (\alpha_0^\omega + \beta_0^\omega)\\
& \label{Qb} \displaystyle
     + \sum_{\omega\in\Omega}w^\omega\sum_{m\in\M^\omega}\sum_{i\in\I_m^\omega\cup\{0\}}\sum_{n\in\N^\omega}\sum_{j\in\J_n^\omega\cup\{0\}}
                            G_{minj}^\omega x_{mi}^\omega y_{nj}^\omega,
\eea
\end{subequations}
\begin{subequations}\label{CDDP-cons}
\bea
& \label{Q1}  \text{s.to }
              \displaystyle \sum_{k\in\K_i}\alpha_{ki} \leq 1         & \forall i\in\I,\\
& \label{Q2}  \displaystyle \sum_{i\in\I}\sum_{k\in\K_i}\alpha_{ki} \leq \overline{I},\\
& \label{Q3}  \alpha_{ki}\in\{0,1\}                                   & \forall k\in\K_i, \, i\in\I,\\
& \label{Q4}  \displaystyle \sum_{k\in\K_j}\beta_{kj} \leq 1          & \forall j\in\J,\\
& \label{Q5}  \displaystyle \sum_{j\in\J}\sum_{k\in\K_j}\beta_{kj} \leq \overline{J},\\
& \label{Q6}  \beta_{kj}\in\{0,1\}                                    & \forall k\in\K_j, \, j\in\J,\\
& \label{Q7}  \displaystyle
              \sum_{m\in\M^\omega:i\in\I_m^\omega}S_m^\omega x_{mi}^\omega \leq (1- D_i^\omega)\sum_{k\in\K_i} S_{ki} \alpha_{ki}
                                                                      & \forall i\in\I, \, \omega\in\Omega,\\
& \label{Q8}  \displaystyle
              x_{m0}^\omega  \leq \alpha_0^\omega,  \,\,
              \sum_{i\in\I_m^\omega\cup\{0\}} x_{mi}^\omega = 1       & \forall m\in\M^\omega, \, \omega\in\Omega,\\
& \label{Q9}  x_{mi}^\omega\in\{0,1\}                                 & \forall i\in\I_m^\omega\cup\{0\}, \, m\in\M^\omega, \, \omega\in\Omega,\\
& \label{Q10} \displaystyle
              \sum_{n\in\N^\omega:j\in\J_n^\omega}R_n^\omega y_{nj}^\omega \leq (1- D_j^\omega)\sum_{k\in\K_j} R_{kj}\beta_{kj}
                                                                      & \forall j\in\J, \, \omega\in\Omega,\\
& \label{Q11} \displaystyle
              y_{n0}^\omega \leq \beta_0^\omega, \,\,
              \sum_{j\in\J_n^\omega\cup\{0\}} y_{nj}^\omega = 1       & \forall n\in\N^\omega, \, \omega\in\Omega,\\
& \label{Q12}  y_{nj}^\omega\in\{0,1\}                                & \forall j\in\J_n^\omega\cup\{0\}, \, n\in\N^\omega, \, \omega\in\Omega,\\
& \label{Q130} \alpha_0^\omega, \, \beta_0^\omega\in\{0,1\}           & \forall \omega\in\Omega.
\eea
\end{subequations}
The objective function \eqref{CDDP-obj} gives the door construction cost in the cross-dock infrastructure
plus the outsourcing penalization plus the expected binary quadratic assignment cost in the scenarios.
The latter includes the standard and outsourcing costs of commodity handling and transportation from strip to stack doors.
The first stage constraints \eqref{Q1}--\eqref{Q6} in system \eqref{CDDP-cons} define the construction of the doors and the installation of the related nominal capacity.
Thus, the clique inequalities \eqref{Q1} and \eqref{Q4} force only one capacity level, if any, for each strip or stack door.
The cover inequalities \eqref{Q2} and \eqref{Q5} bound the number of doors to be built.
The binary domain of the variables is defined in \eqref{Q3} and \eqref{Q6}.
The second stage constraint system \eqref{Q7}--\eqref{Q130} defines the assignment constraints in the cross-dock for each scenario.
Thus, \eqref{Q7} and \eqref{Q10} upper-bound the commodity volume handled by the strip and stack doors, respectively, based on the doors' net capacity.

Given the difficulty of solving the combinatorial problem CDDP-TS and the frequent requirement for the high cardinality of the uncertain sets,
a mathematically equivalent Linearization mixed Integer Programming (LIP) is proposed, and additionally,
a scenario Cluster Decomposition (CD) of model LIP  into smaller submodels is also proposed.
The aim is to obtain (hopefully, strong) lower and upper bounds of the CDDP-TS solution value (see the next section).

\section{LIP reformulation for the binary quadratic model CDDP-TS} \label{sec:LIP}
Model LIP requires to replace the quadratic term $x_{mi}^\omega y_{nj}^\omega$ in expression \eqref{Qb} with the continuous variable $v_{minj}$,
for $m\in\M^\omega, \, i\in\I^\omega\cup\{0\}, \, n\in\N^\omega, \, j\in\J^\omega\cup\{0\}, \, \omega\in\Omega$,
by considering RLT1; see \cite{SheraliAdams94,Guignard12} and others.
Now, since it is required that $x_{mi}^\omega\cdot y_{nj}^\omega \, = \, v_{minj}^\omega$
(and thus, $x_{mi}^\omega=0 \vee y_{nj}^\omega=0 \,\Longleftrightarrow \, v_{minj}^\omega=0$ and $x_{mi}^\omega=1 \wedge y_{nj}^\omega=1 \, \Longleftrightarrow \, v_{minj}^\omega=1$),
both sides of the constraints \eqref{Q8} and \eqref{Q11} can be multiplied by $y_{nj}$ and $x_{mi}$, respectively.
Thus, the new constraints \eqref{LIP2}--\eqref{LIP4} are appended in order to satisfy the requirements for equivalency.

Model LIP can be expressed
\begin{equation}\label{LIP-obj}
 z^* \, = \, \min \quad \eqref{Qa} +
     \displaystyle \sum_{\omega\in\Omega}w^\omega\sum_{m\in\M^\omega}\sum_{i\in\I_m^\omega\cup\{0\}}\sum_{n\in\N^\omega}\sum_{j\in\J_n^\omega\cup\{0\}}
                   G_{minj}^\omega v_{minj}^\omega,
\end{equation}
\begin{subequations}\label{LIP-cons}
\bea
& \label{LIP1} \text{s.to cons system } \eqref{CDDP-cons},\\
& \label{LIP2} \displaystyle \sum_{j\in\J_n^\omega\cup\{0\}} v_{minj}^\omega = x_{mi}^\omega
                     & \forall i\in\I_m^\omega\cup\{0\}, \, m\in\M^\omega,  n\in\N^\omega, \, \omega\in\Omega,\\
& \label{LIP3} \displaystyle \sum_{i\in\I_m^\omega\cup\{0\}} v_{minj}^\omega = y_{nj}^\omega
                     & \forall m\in\M^\omega, \, j\in\J_n^\omega\cup\{0\}, \, n\in\N^\omega, \, \omega\in\Omega,\\
& \label{LIP4} v_{minj}^\omega\in[0,1]
                     & \forall i\in\I_m^\omega\cup\{0\}, \, m\in\M^\omega,  j\in\J_n^\omega\cup\{0\}, \, n\in\N^\omega, \, \omega\in\Omega.
\eea
\end{subequations}
The objective function \eqref{LIP-obj} is the same as \eqref{CDDP-obj},
where the binary quadratic terms \eqref{Qb} are replaced with continuous linear terms.
The constraint system \eqref{LIP-cons} is composed of system \eqref{CDDP-cons} and the new one \eqref{LIP2}--\eqref{LIP4} related to the $(x,y,v)$-variables.

\section{Scenario Cluster Decomposition of model LIP}\label{sec:cd}
There are many decomposition algorithms for solving the large-scale instances of general two-stage mixed-integer linear models;
see in \cite{mcldtsd17} an overview of exact and matheuristic algorithmic typologies in the literature for problem solving.
Moreover, there is a high computational complexity in the special structure of model LIP, as well as in the structure of the deterministic CDAP model; see \cite{cdap23}.
The latter is embedded in the CDDP-TS constraint system for each scenario, once the first stage variables are fixed.
Thus, the development of an ad-hoc decomposition approach is motivated, and this section is organized as follows:
Subsection \ref{sec:cluster-gen} introduces the approach for scenario cluster generation.
Subsection \ref{sec:c-submodels} presents two options for the scenario cluster $c$-submodel that provides a lower bound of the solution value of model LIP, given a set of scenario clusters represented by $\C$.
Subsection \ref{sec:lazy} presents a lazy algorithm for providing a feasible solution of model LIP, given  first stage feasible $\alpha$- and $\beta$-vectors.

\subsection{Scenario cluster generation}\label{sec:cluster-gen}
A scenario Cluster Decomposition (CD) scheme is considered to provide lower bounds of the solution value of model LIP, as well as a basis for obtaining feasible solutions for the model.
The approach consists of reformulating the model by taking advantage of the fact that the uncertainty also depends on the sets and not only on the parameters.
The two-step scheme that is proposed for scenario cluster generation is as follows:
\begin{enumerate}
\item Let $\Omega^{c1} \subseteq \Omega$ denote the set of scenarios in cluster $c1$, for $c1\in\C^1$, where
$\C^1$ is the set of clusters, such that
$\Omega = \cup_{c1\in\C^1} \Omega^{c1}, \,\, \Omega^{c1}\cap\Omega^{c1'} =\emptyset, c1,c1'\in\C^1:c1 \neq c1'$,
where two scenarios $\omega$ and $\omega'$ are assigned to the same cluster $c1$ provided they have the same number of origin and destination nodes, i.e. $\Omega^{c1} =\{\omega,\omega'\in\Omega: |\M^\omega|=|\M^{\omega'}| \, \wedge |\N^\omega|=|\N^{\omega'}|\}$.

\item If $|\Omega^{c1}|>1$, i.e. cluster $c1$ is non-singleton, for $c1\in\C^1$, it can be
partitioned into the cluster subset $\C_{c1}^2$, such that  $\Omega^{c1} = \cup_{c2\in\C_{c1}^2}\Omega_{c1}^{c2}$ and $|\Omega_{c1}^{c2}| \leq \kappa_{c1}$, where $\kappa_{c1}$ is a modeler-driven parameter that gives the maximum number of scenarios in any cluster $c2\in\C_{c1}^2$.
As a matter of fact, $|\Omega_{c1}^{c2}| = \kappa_{c1}$, for all clusters $c2$ but cluster $\ell$, where $|\Omega_{c1}^\ell| < \kappa_{c1}$,
such that $\Omega_{c1}^\ell = \Omega^{c1} \setminus \bigcup_{c2\in\C_{c1}^2: c2 \neq \ell} \Omega_{c1}^{c2}$.

\item The assignment of the scenarios to set $\Omega_{c1}^{c2}$ is performed \textit{at random} from set $\Omega^{c1}$.
\end{enumerate}
For simplification purposes, let us group all clusters into one set $\C$, so that $\C=\{c2\in\C_{c1}^2 \, \forall c1\in\C\}$.

\subsection{Lower bound-related submodels}\label{sec:c-submodels}
Two relaxation submodels are considered to obtain a lower bound of model LIP \eqref{LIP-obj}--\eqref{LIP-cons} for a given scenario cluster set $\C$. Option 1 consists of a CD of the first stage constraints and variables.

The notation for the new elements to be required in CD is as follows, for $c\in\C$:
\begin{description}\parskip 0.5mm
\item[$\I^c$,] subset $\{i\}$ of strip doors expressed as
     $\{ i\in\I: \exists m\in\M^\omega, \text{ where } i\in\I_m^\omega, \, \omega\in\Omega^c \}$.
\item[$\alpha_{ki}^c$,] copy of variable $\alpha_{ki}$, for $k\in\K_i, \, i\in\I^c$.
\item[$\J^c$,] subset $\{j\}$ of stack doors expressed as
     $\{ j\in\J: \exists n\in\N^\omega, \text{ where } j\in\J_n^\omega, \, \omega\in\Omega^c \}$.
\item[$\beta_{kj}^c$,] copy of variable $\beta_{kj}$, for $k\in\K_j, \, j\in\J^c$.
\item[$w^c$,] weight factor representing the likelihood that is associated with scenario cluster $c$, for $c\in\C$, such that
    $w^c = \sum_{\omega\in\Omega^c}w^\omega$; then, $\sum_{c\in\C}w^c=1$.
\item[$w'^\omega$,] weight factor representing the likelihood that is associated with scenario $\omega$ in cluster $c$,
    for $\omega\in\Omega^c, \, c\in\C$, such that
    $w'^\omega = \frac{w^\omega}{w^c}$; then, $\sum_{\omega\in\Omega^c}w'^\omega=1$.
\end{description}

\subsubsection{Option 1. First-stage CD $c$-submodel of model LIP, for $c\in\C$}
Let $\underline{z}^c$ be the solution value of the CD $c$-submodel.
It can be expressed
\begin{subequations}\label{c-model}
\bea
& \label{lb-cost} \displaystyle
    \underline{z}^c\, = \, \min \sum_{i\in\I^c}\sum_{k\in\K_i} F_{ki} \alpha_{ki}^c + \sum_{j\in\J^c}\sum_{k\in\K_j}F_{kj} \beta_{kj}^c
 +  \sum_{\omega\in\Omega^c}w'^\omega F_0 (\alpha_0^\omega + \beta_0^\omega)\\
& \label{lb-cost-1} \displaystyle
 +  \sum_{\omega\in\Omega^c}w'^\omega\sum_{m\in\M^\omega}\sum_{i\in\I_m^\omega\cup\{0\}}\sum_{n\in\N^\omega}\sum_{j\in\J_n^\omega\cup\{0\}}
                            G_{minj}^\omega v_{minj}^\omega,\\
& \text{s.to cons system } \eqref{LIP-cons}, \text{ where } \alpha_{ki}, \beta_{kj}, \text{ and } \Omega
                  \text{ are replaced with } \alpha_{ki}^c, \beta_{kj}^c, \text{ and } \Omega^c, \text{ resp.}
\eea
\end{subequations}
Thus, a lower bound of model LIP can be expressed
\begin{equation}\label{lb}
\displaystyle \underline{z} = \sum_{c\in\C}w^c \underline{z}^c.
\end{equation}

\subsubsection{Option 2. First-stage CD $c$ strip- and stack-based scenario $c$-submodels, for $c\in\C$}\label{sec:c-strip-stack}
This subsection presents a deeper relaxation of model LIP,
where the stack- and strip-door submodels are independently considered;
let us call them CD strip $c$-submodel \eqref{strip-c-model} and CD stack $c$-submodel \eqref{stack-c-model},
such that $\underline{z}_{strip}^c$ and $\underline{z}_{stack}^c$ are the related solution values, for $c\in\C$:

{\it CD strip $c$-submodel}
\begin{subequations}\label{strip-c-model}
\bea
& \label{lb-strip-obj} \underline{z}_{strip}^c \, = \, \min \, \displaystyle
     \sum_{i\in\I^c}\sum_{k\in\K_i}F_{ki} \alpha_{ki}^c
   + \sum_{\omega\in\Omega^c}w'^\omega F_0 \alpha_0^\omega\\
& \label{lb-strip-obj-b} \displaystyle
   + \sum_{\omega\in\Omega^c}w'^\omega\sum_{m\in\M^\omega}\sum_{i\in\I_m^\omega\cup\{0\}}\sum_{n\in\N^\omega}\sum_{j\in\J_n^\omega\cup\{0\}}
             \frac{1}{2}G_{minj}^\omega v_{minj}^\omega,\\
& \nonumber
\text{s.to cons system } \eqref{Q1}-\eqref{Q3}, \eqref{Q7}-\eqref{Q9}, \eqref{LIP2}, \eqref{LIP4}, \text{ where } \\
& \label{lb-strip-Q} \qquad \alpha_{ki}, \, \I, \, \Omega
                     \text{ are replaced with } \alpha_{ki}^c, \,  \I^c, \, \Omega^c,\\
& \label{lb-strip-Q130} \alpha_0^\omega\in\{0,1\}  \quad \forall m\in\M^\omega, \, \omega\in\Omega^c.
\eea
\end{subequations}

{\it CD stack $c$-submodel}
\begin{subequations}\label{stack-c-model}
\bea
& \label{lb-stack-obj} \underline{z}_{stack}^c \, = \, \min \, \displaystyle
     \sum_{j\in\J^c}\sum_{k\in\K_j} F_{kj} \beta_{kj}^c
   + \sum_{\omega\in\Omega^c}w'^\omega F_0 \beta_0^\omega\\
& \label{lb-stack-obj-b} \displaystyle
   + \sum_{\omega\in\Omega^c}w'^\omega\sum_{m\in\M^\omega}\sum_{i\in\I_m^\omega\cup\{0\}}\sum_{n\in\N^\omega}\sum_{j\in\J_n^\omega\cup\{0\}}
             \frac{1}{2}G_{minj}^\omega {v'}_{minj}^\omega,\\
& \nonumber
\text{s.to cons system } \eqref{Q4}-\eqref{Q6}, \eqref{Q10}-\eqref{Q12}, \eqref{LIP3}, \eqref{LIP4}, \text{ where } \\
& \label{lb-stackQ} \qquad \beta_{kj}, \, \J, \, \Omega \, \text{ are replaced with }
                    \beta_{kj}^c, \, \J^c, \, \Omega^c,\\
& \label{lb-stack-Q130} \beta_0^\omega\in\{0,1\}  \quad \forall n\in\N^\omega, \, \omega\in\Omega^c.
\eea
\end{subequations}
Thus, another lower bound of model LIP can be expressed
\begin{equation}\label{lb-week}
\displaystyle \underline{\underline{z}} = \sum_{c\in\C}w^c( \underline{z}_{strip}^c + \underline{z}_{stack}^c ).
\end{equation}

\subsection{Upper bound-related submodel}\label{sec:lazy}
A lazy heuristic (LH) scheme is proposed for obtaining feasible solutions of the model CDDP-TS \eqref{CDDP-obj}--\eqref{CDDP-cons},
based on an $(\alpha, \, \beta)$-solution that is assumed to satisfy the first stage constraints \eqref{Q1}--\eqref{Q6}.
The procedure basically consists of fixing the $\alpha$- and $\beta$-variables to the $\hat{\alpha}$- and $\hat{\beta}$-solutions, respectively, that can be  retrieved by solving the $c$-submodels.

Consider the following notation:
\begin{description}\parskip 0.5mm
\item[$k(i)$] and $k(j)$, capacity levels in strip door $i$ and stack door $j$ that belong to the solutions $\alpha_{ki}=1$ and $\beta_{kj}=1$ in the $c$-submodels, for  $k\in\K_i, \, i\in\I$ and $k\in\K_j, \, j\in\J$, respectively.
Note that by construction, these values satisfy the first stage constraints \eqref{Q1}--\eqref{Q6} in system \eqref{CDDP-cons}.
\end{description}

The notation of the fixed $\hat{\alpha}$ and $\hat{\beta}$ of the $\alpha$- and $\beta$-variables to be used in the scenario $\omega$-submodel \eqref{h-omega} is as follows:
\begin{description}\parskip 0.5mm
\item[$\hat{\alpha}_{k(i)i}=1$] and $\hat{\beta}_{k(j)j}=1$, for $i\in\I$ and $j\in\J$, respectively, and

\item[$\hat{\alpha}_{k(i)i}=1$,] where $k(i)=0$, for $i\in\I: \, \alpha_{ki}=0 \, \forall \, k\in\K_i$, and
      $\hat{\beta}_{k(j)j}=1$,  where $k(j)=0$,  for $j\in\J: \, \beta_{kj}=0  \, \forall \, k\in\K_j$,
      so that the \textit{basic capacities} $S_{0i}$ and $R_{0j}$ are considered in the scenario $\omega$-submodels \eqref{h-omega} to help to avoid outsourcing solutions.
\end{description}

As a result, the second stage submodel of the original one is trivially decomposed into a set of $|\Omega|$ independent scenario $\omega$-submodels to be expressed
\begin{subequations}\label{h-omega}
\bea
& \label{h1} \displaystyle z_2^\omega \, = \, \min \quad F_0 (\alpha_0^\omega + \beta_0^\omega)
     + \sum_{m\in\M^\omega}\sum_{i\in\I_m^\omega\cup\{0\}}\sum_{n\in\N^\omega}\sum_{j\in\J_n^\omega\cup\{0\}}
            G_{minj}^\omega x_{mi}^\omega y_{nj}^\omega,\\
& \text{s.to cons system } \eqref{Q8}, \eqref{Q9}, \eqref{Q11}-\eqref{Q130}, \text{ but only for the given scenario } \omega,\\
& \label{h7}   \displaystyle
               \sum_{m\in\M^\omega:i\in\I_m^\omega}S_m^\omega x_{mi}^\omega \leq (1- D_i^\omega) S_{k(i)i}
                                                                             \quad \forall i\in\I,\\
& \label{h10}  \displaystyle \sum_{n\in\N^\omega:j\in\J_n^\omega}R_n^\omega y_{nj}^\omega \leq (1- D_j^\omega) R_{k(j)j}
                                                                             \quad \forall j\in\J.
\eea
\end{subequations}
Note that the submodel is the CDAP binary quadratic $\omega$-model, for $\omega\in\Omega$.
Thus, the local search heuristic (LSH) considered in \cite{cdap23} can be used to obtain a (hopefully, good) feasible solution for large-scale instances;
otherwise, a state-of-the-art solver can be used for either the binary quadratic model above or its LIP reformulation.

\section{Algorithm SCS4B, scenario clustering and splitting for a bounding solution of CDDP-TS}\label{sec:SCS4B}
A matheuristic algorithm is introduced for obtaining feasible solutions and the optimality gap of the incumbent solution value for model CDDP-TS \eqref{CDDP-obj}--\eqref{CDDP-cons}.
It has an aversion to solutions in which outsourcing is required (i.e. solutions with $\alpha_0^\omega=1$ or $\beta_0^\omega=1$, for any scenario $\omega\in\Omega$).

Roughly, Algorithm \ref{alg} has the following main items:
The first is devoted to refining the scenario set $\Omega$ and generating the singleton cluster set $\S$ (Step 0).
For this purpose, the CD submodel presented in Subsection \ref{sec:c-submodels} is solved, where $c=\{\omega\}$, for $\omega\in\Omega$.
If $\exists \omega\in\Omega$: $\hat{x}_0^\omega=1 \vee \hat{y}_0^\omega=1$
(i.e. it is an outsourcing solution and thus a non-feasible one for our purposes),
then scenario $\omega$ is removed from set $\Omega$.
Thus, $\S$ is a set with only non-outsourcing scenarios, and the related first stage $(\hat{\alpha}, \, \hat{\beta})$-solutions are saved in pool $\S$.
In Step 1, the scenario cluster set $\C$ is generated.
If the outsourcing solution is retrieved from the $c$-submodel for the general cluster set $\C$ in Step 2,
then cluster $c$ is partitioned into singleton ones and thus the outsourcing scenarios are removed from the sets $\C$ and $\S$.

Step 3 performs the following tasks:
\begin{itemize}
\item Fixing, in model CDDP-TS, the first stage $(\alpha, \, \beta)$-variables to the $(\hat{\alpha}, \, \hat{\beta})$-values retrieved from the solution obtained in the $c$-submodels, for $c\in\C$, plus the \textit{basic capacities} for the door sets  $\{i\in\I^c: \, k(i)=0\}$ and $\{j\in\J^c: \, k(j)=0\}$;

\item Solving the second stage-related $\omega$-submodels \eqref{h-omega}, for $\omega\in\Omega$,
which solutions of the $(x, \, y)$-variables are retrieved from;  and

\item Computing the related solution value.
For this purpose, let the new door sets $\hat{\I}^c$ and $\hat{\J}^c$ denote the strip and stack doors that have been used in, at least, one scenario $\omega$-submodel.  They are expressed as
$\hat{\I}^c = \{i\in\I^c: \, \exists m\in\M^\omega, \, \omega\in\Omega: \, \hat{x}_{mi}=1\}$ and
$\hat{\J}^c = \{j\in\J^c: \, \exists n\in\N^\omega, \, \omega\in\Omega: \, \hat{y}_{nj}=1\}$,
such that $(\hat{x}, \, \hat{y})$ is the solution vector, and the cost $z_1^c$ of the installation of the door infrastructure, as well as the total cost, can be expressed as follows:
\bea
& \label{z_1} \displaystyle z_1^c = \sum_{i\in\hat{\I}^c}F_{k(i)i} \hat{\alpha}_{k(i)i} + \sum_{j\in\hat{\J}^c}F_{k(j)j} \hat{\beta}_{k(j)j},\\
& \label{z-c} \displaystyle z^c=z_1^c+\sum_{\omega\in \Omega}w^\omega z_2^\omega.
\eea
\end{itemize}

Step 4 evaluates the goodness of the new solution for original model CDDP-TS, as well as its incumbency.
If the solution has outsourcing or if it is not the incumbent one, then,
the capacity of the doors $\{i\in\I^c: \, k(i)>0\}$ and $\{j\in\J^c: \, k(j)>0\}$ is increased for the current cluster $c$;
afterwards, Step 3 is executed again for that cluster.
Otherwise, the incumbent solution is updated.
Step 5 determines if there is still a cluster that has not yet been explored, in which case the procedure goes to Step 3;
otherwise, the procedure is over.

\begin{algorithm*}
  \begin{footnotesize}
  \caption{SCS4B procedure for obtaining feasible solutions and the optimality gap in model CDDP-TS \eqref{CDDP-obj}--\eqref{CDDP-cons}}\label{alg}
\vskip 0.2cm
  \begin{description}
\item [{\bf Step 0:}]\textit{Aim: Generate the pool $\S$ of singleton scenario clusters, where $c=\{\omega\}$, for $\omega\in\Omega$.}
\begin{itemize}\parskip 0.5mm
\item Solve the modeler-driven choice between $c$-submodel \eqref{c-model} and the pair of strip and stack $c$-submodels
\eqref{strip-c-model} and \eqref{stack-c-model}.
    \\Note: For simplification purposes, let us assume that $c$-submodel \eqref{c-model} is the modeler's submodel.
\item If the solution is an outsourcing one, then remove scenario $\omega$ from set $\Omega$ and pool $\S$.
    Otherwise, update the pool with the first stage solution $(\hat{\alpha}_{k(i)i} \, \forall i\in\I^c, \, \hat{\beta}_{k(j)j} \, \forall j\in\J^c)$
    retrieved from the $c$-submodel.
\end{itemize}

\item [{\bf Step 1:}]\textit{Aim: Generate the scenario cluster set $\C$}, as presented in Subsection \ref{sec:cluster-gen}.

\item [{\bf Step 2:}] \textit{Aim: Obtain the first stage $(\hat{\alpha}, \,\hat{\beta})$-solution for the scenario clusters in set $\C$.}
\\Select a not-yet-explored cluster $c\in\C$.
\begin{itemize}\parskip 0.5mm
\item If $|\Omega^c|=1$ (i.e. cluster $c$ is singleton), retrieve the first stage solution
    $\hat{\alpha}_{k(i)i} \, \forall i\in\I^c$ and $\hat{\beta}_{k(j)j} \, \forall j\in\J^c$ from pool $S$.

\item Otherwise, retrieve those values from the solution of $c$-submodel \eqref{c-model},
    and if it is a non-outsourcing one (i.e. $\hat{\alpha}_0^\omega=0$ and $\beta_0^\omega=0$ $\forall \omega\in\Omega^c$),
    then set up
    $\hat{\alpha}_{k(i)i}:=\alpha_{k(i)i}^{c} \, \forall i\in\I^c$ and
    $\hat{\beta}_{k(j)j}:={\beta}_{k(j)j}^{c}  \, \forall j\in\J^c$.

\item If the solution is an outsourcing one, then cluster $c$ is partitioned into singleton subclusters.
    Thus, $c$ is removed from set $\C$ and the scenarios in $\Omega^c$ are appended to set $\C$ under the label `already explored',
    since the $(\hat{\alpha}, \, \hat{\beta})$-solution is taken from the singleton's pool $\S$.
\end{itemize}

\item [{\bf Step 3:}] \textit{Aim: Obtain the second stage $(\hat{x}, \, \hat{y})$-solution for the scenario clusters in set $\C$.}
\begin{itemize}\parskip 0.5mm
\item For a given cluster $c\in\C$, append the fixed first stage variables $\hat{\alpha}:=\alpha^c, \, \hat{\beta}:=\beta^c$ in model CDDP-TS.

\item Thus, it is trivially decomposed into a set of $|\Omega|$ independent second stage scenario submodels,
    which related solution value is $z^c$ \eqref{z-c} (see Subsection \ref{sec:lazy}).
\end{itemize}

\item [{\bf Step 4:}] \textit{Aim: Evaluate the goodness of the first stage solution $( \hat{\alpha}_{k(i)i}, \, \hat{\beta}_{k(j)j} )$ from cluster $c\in\C$.}
\begin{itemize}\parskip 0.5mm
\item Compute $out$ as $\sum_{\omega\in\Omega}(\hat{\alpha}_0^\omega + \hat{\beta}_0^\omega)$, where the solutions are retrieved from the $\omega$-submodels \eqref{h-omega} $\forall \omega\in\Omega$.
\item If $out \leq \rho |\Omega|$ and $ \hat{z}^c < \overline{z}$, then
    update the incumbent solution value $\overline{z}$, as well as the first- and second stage solution vectors:
    \begin{description}\parskip 0.5mm
    \item $\overline{z}:= \hat{z}^c$;
        $\overline{\alpha}_{k(i)i}:= \hat{\alpha}_{k(i)i}, \, \overline{\beta}_{k(j)j}=\hat{\beta}_{k(j)j}$;\\
        $\overline{x}_{mi}^\omega:= \hat{x}_{mi}^\omega$ for $i\in\I_m^\omega\cup\{0\}, \, m\in\M^\omega, \, \omega\in\Omega$;\\
        $\overline{y}_{nj}^\omega:= \hat{y}_{nj}^\omega$ for $j\in\J_n^\omega\cup\{0\}, \, n\in\N^\omega, \, \omega\in\Omega$; and
       go to \textbf{Step 5}.
       \end{description}
\item Otherwise, if $\exists i\in\I^c: \, 0<k(i) \neq |\K_i|$ or $\exists j\in\J^c: 0<k(j) \neq |\K_j|$, then
      fix the $\alpha$- and $\beta$-variables at 1 for higher capacity levels.
      Thus, $k(i):=min\{ k(i) + \delta \cdot |\K_i|, \, |\K_i| \}$,
      $k(j):=min\{ k(j) + \delta \cdot |\K_j|, \, |\K_j| \}$,
      $\hat{\alpha}_{k(i)i}:=1$, for $i\in\I^c$, and $\hat{\beta}_{k(j)j}:=1$, for $j\in\J^c$, and go to \textbf{Step 3}.
\end{itemize}

\item [{\bf Step 5:}]\textit{Check if all scenario clusters have been explored.}
\begin{itemize}\parskip 0.5mm
\item If there is a scenario cluster in set $\C$ for which first stage solution fixing has not yet been explored,
then go to \textbf{Step 3}.
\item Otherwise, the optimality gap is $100.\frac{\overline{z}-\underline{z}}{\overline{z}}$, and the algorithm stops,
where there is a choice between $\underline{z}$ \eqref{lb} and $\underline{z}=\underline{\underline{z}}$ \eqref{lb-week}.
\end{itemize}
\end{description}
  \end{footnotesize}
\end{algorithm*}

Consider the following additional notation required by Algorithm \ref{alg}:
\begin{description}\parskip 0.5mm
\item[$out$,] number of scenarios with `outsourcing' in set $\Omega$, i.e. the solutions that have either $\hat{\alpha}_0^\omega=1$ or $\hat{\beta}_0^\omega=1$, for given fixed values of the $\alpha$- and $\beta$-variables.
\item[$\rho$,] modeler-driven maximum fraction of scenarios with `outsourcing' in set $\Omega$ to consider that the solution is still a feasible one.
\item[$\overline{z}$,] incumbent solution value of model CDDP-TS \eqref{CDDP-obj}--\eqref{CDDP-cons}.
    Let its initialization be $\overline{z} = \infty$.
\item[$\delta \geq 1$,] number of capacity levels to increase the current capacity one for $k(i)>0, \, i\in\I$ and $k(j)>0, \, j\in\J$ if a non-feasible solution is found for set $\Omega$ in a given iteration.
\end{description}

\section{Computational experiments}\label{sec:results}
The experiments reported in this section were conducted using a
Debian Linux workstation (kernel v4.19.0), ES-2670 processors (20 threads, 2.50 Ghz), and 128 GB RAM.
The CPLEX v12.8 Concert Technology library has been used, embedded in
C/C++  experimental code, to solve the decomposition submodels in the CD and LH schemes, as well as to straightforwardly solve the full models.

\subsection{Testbed generation}\label{sec:testbed-gen}
A testbed of instances is generated by considering the stochastic extension of some of the deterministic instances
of model CDAP that were considered in \cite{cdap23}.
This subsection presents the related scheme.
As a matter of fact, most of the deterministic CDAP instances that have been considered to generate the stochastic testbed were proposed in \cite{Guignard12}.
They have a similar structure, which makes it possible to embed them as scenarios for the second stage in model CDDP-TS (see below),
given their common framework of parameters and variables for the first stage.
The dimensions of these instances are
$\{8,9,10,11,12,20,25\}$ for $|\M|=|\N|$ and $\{4,5,6,7,10\}$ for $|\I|=|\J|$.
As in \cite{Guignard12}, the nonzero (25\%) elements of the flow square matrix
are \textit{randomly} generated in the range from 10 to 50.
On the other hand, the elements in the unit distance matrix are in the
range from $8$ to $8+|\I|-1$, where a symmetric pattern is considered for the distance between strip door $i$ and stack door $j$.
Additionally, the capacity of the doors in each instance is set to the total flow though all origins and destinations
divided by the total number of strip and stack doors,
plus a slackness percentage $S$ of the total flow, where $S\in\{5,10,15,20,30\}$.

The so-named \textit{Basic Scenario Cluster} (BSC) instance $|\M|\times|\I|S$ can be generated as follows:
\begin{itemize}\parskip 0.5mm
\item The first stage constraints and variables and the ones related to the inbound and outbound flows for the first scenario in the second stage are taken from the related deterministic CDAP model $|\M|\times|\I|Ss'$ in \cite{cdap23},
    where $s'$ denotes the slackness percentage of the total flow.

\item The other scenarios in the second stage are built by replacing $s'$ with $s$, for $s\in\{5,10,15,20,30\}\setminus\{s'\}$.

\item Additionally, the `outsourcing' doors are penalized in the objective function.
\end{itemize}

\begin{figure}[h]
\begin{center}
\includegraphics[width=10cm, viewport=5 470 315 644,clip=]{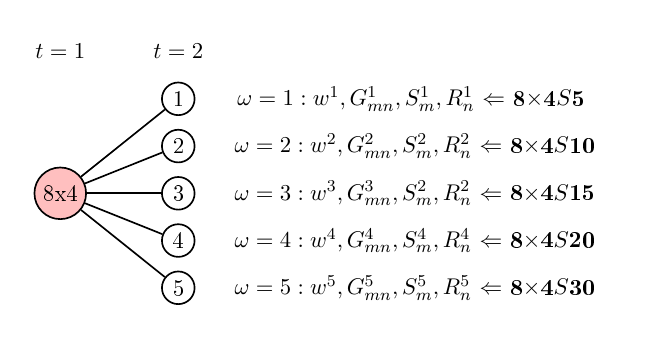}
\vskip -0.9cm
\caption{BSC instance $8\times4S$} \label{fig-cluster}
\end{center}
\end{figure}

The testbed is composed of some of the BSC instances in the set $\{8\times4S, 10\times5S, 15\times6S, 20\times10S\}$ plus the larger instances that result from their total or partial merging.
For illustrative purposes, Fig. \ref{fig-cluster} depicts the smallest BSC instance $8\times 4S$ in the experiment,
where $|\I|=|\J|=4$, $|\M^\omega|=|\N^\omega|=8$ $\forall \omega\in\Omega=\{1,...,5\}$.

It is worth noting that a merging of BSC instances requires some type of homogenization for the sets of strip and stack doors inherited in the enlarged BSC instance obtained from the selected BSC instances.
Consider the following notation:

\begin{description}\parskip 0.5mm
\item[$\BB$,] set $\{b\}$ of the selected BSC instances.

\item[$\I_b, \, \J_b$,] sets of strip and stack doors in BSC $b$, respectively, for $b\in\BB$.

\item[$\Omega_b$,] set of scenarios inherited from BSC $b$, for $b\in\BB$.

\end{description}
The door capacity disruption fractions are fixed as follows for homogenisation purposes:
$D_i^\omega=1$, for $i\in \I\setminus\I_{b}$ (respectively, $D_j^\omega=1$, for $j\in \J\setminus\J_{b}$), for $\omega\in\Omega\setminus\Omega_b,b\in\BB$.

\begin{figure}[h]
\begin{center}
\includegraphics[width=12cm, viewport=16 569 443 786,clip=]{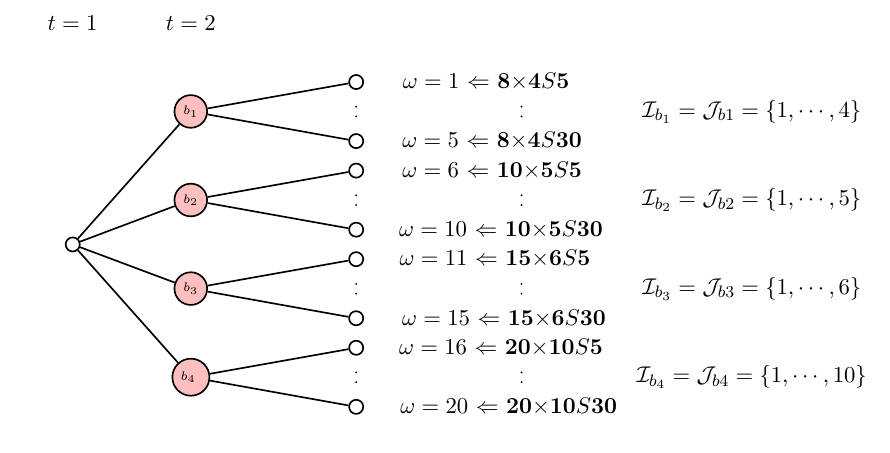}
\vskip -0.8cm
\caption{General structure of the scenario tree of instance I7, built from the integration of BSC nodes} \label{fig-cluster-d}
\end{center}
\end{figure}

Additionally, the first stage parameters (i.e. the nominal capacity for each level $k$ of the strip and stack doors) require,
for feasibility reasons, some adjustments of their values in order to
embed basic BSC instances into a larger instance; this ensures that no outsourcing situation occurs.
Fig. \ref{fig-cluster-d} depicts the two-stage scenario tree based
on the integration of some BSCs that extend some deterministic instances in a stochastic framework.

The headings of Table \ref{tb:1}, where the problem's elements and dimensions are shown, are as follows:
$Inst$, code; $|\Omega|$, number of scenarios;
$|\I|$, number of strip doors; $|\J|$, number of stack doors; and
$|\underline{\M}^\omega|$ and $|\overline{\M}^\omega|$ (resp., $|\underline{\N}^\omega|$ and $|\overline{\N}^\omega|$),
minimum and  maximum number of origins (respectively,  destinations) in the scenarios.
The following fixed values are considered in the experiment for all the instances:
$|K_i|=|\K_j|=5,\, i\in\I,\, j\in\J$; $\overline{I}= |\I|+1$; and $\overline{J}= |\J|+1$.

\begin{table}[!h]
\small
  \caption{\footnotesize{Problem CDDP-TS. Dimensions}}
  \label{tb:1}
 \vskip .2cm
\begin{tabular}{l|rrrrr}
  \hline
 &&&&&\\[-1.5ex]
$Inst$ & $|\Omega|$ & $|\I|$ & $|\J|$ & $|\underline{\M}^\omega|$-$|\overline{\M}^\omega|$ & $|\underline{\N}^\omega|$-$|\overline{\N}^\omega|$\\[0.5ex]
\hline
I1:$\,\, 8\times 4S$
&5&4&4&8-8&8-8\\
I2:$\,\, 10\times 5S$
&5&5&5&10-10&10-10\\
I3:$\,\, 8\times 4S\cup 10\times 5S$
&10&5&5&8-10&8-10\\
I4:$\,\, 15\times 6S$
&5&6&6&15-15&15-15\\
I5:$\,\, 8\times 4S\cup 10\times 5S\cup 15\times 6S$
&15&6&6&8-15&8-15\\
I6:$\,\, 20\times 10S$
&5&10&10&20-20&20-20\\
I7:$\,\, 8\times 4S\cup 10\times 5S\cup 15\times 6S \cup 20\times 10S$
&20&10&10&8-20&8-20\\
I8:$\,\, 8\times 4S\cup 10\times 5S\cup 15\times 6S \cup 20\times 10S$
&60&10&10&8-20&8-20\\
\hline
\end{tabular}
\end{table}

\subsection{Straightforward CPLEX results}
The dimensions of the stochastic two-stage MILP model (model LIP) \eqref{LIP-obj}--\eqref{LIP-cons} are shown in Table \ref{tb:2}, which has the following headings:
$\#m$, number of constraints;
$\#n01$, number of binary variables;
$\#nc$,  number of continuous variables;
$\#nz$, number of nonzero elements;
$|\C|$, number of clusters; and
$|\Omega^c|$, number of scenarios in the largest cluster.
Note that the subindex $c$ refers to the dimensions of the largest
scenario cluster considered, for $c\in \C$.

\begin{table}[!h]
\small
  \caption{\footnotesize{Model LIP \eqref{LIP-obj}--\eqref{LIP-cons}. Dimensions}}
  \label{tb:2}
 \vskip .2cm
\begin{tabular}{l|rrrr|rr|rrrr}
  \hline
 &&&&&&&&&&\\[-1.5ex]
$Inst$ & $\#m$ & $\#n01$ & $\#nc$ & $\#nz$ & $|\C|$&$|\Omega^c|$ & $\#m_c$ & $\#n01_c$ & $\#nc_c$ & $\#z_c$\\[0.5ex]
\hline
 I1&3410&450&8000&20360&2&3&2050&286&4800&12248\\
 I2&6262&660&18000&43650&3&2&2512&294&7200&17520\\
 I3&9662&1070&26000&63930&5&2&2512&294&7200&17520\\
 I4&16124&1120&55125&128670&5&1&3236&272&11025&25830\\
 I5&25774&2140&81125&192500&9&5&3414&470&8000&20400\\
 I6&44522&2310&242000&533300&2&4&35622&1868&193600&426680\\
  I7&70282&4390&323125&725480&3&4&35662&1868&193600&426680\\
  I8&210802&12970&969375&2176640&3&10&89022&4520&484000&1066400\\
\hline
\end{tabular}
\end{table}

Table \ref{tb:3} provides the main results obtained when model LIP \eqref{LIP-obj}--\eqref{LIP-cons} is solved by straightforward CPLEX.
The headings are as follows: $z_{LP}$ and $t_{LP}$, LP lower bound and wall time (seconds);
$\underline{z}_{CPX}$, lower bound of the optimal solution value;
$\overline{z}_{CPX}$ and $t_{\overline{z}_{CPX}}$, incumbent solution value and wall time; and
$GAP_{CPX}\%$, optimality gap of the incumbent solution, expressed as
$100\cdot {{\overline{z}_{CPX}-\underline{z}_{CPX}}\over {\overline{z}_{CPX}}}$.
Note that CPLEX reaches the time limit of 43200 seconds (12 hours) in all instances but the smallest one, I1, without guaranteeing optimality.
It is worth noting that the time reported in the brackets, as shown in the column under the heading $t_{\overline{z}_{CPX}}$,
is the amount of time that the CPLEX incumbent solution required to improve the solution value obtained by the matheuristic SCS4B.

\begin{table}[!h]
\small
  \caption{\footnotesize{Model LIP \eqref{LIP-obj}--\eqref{LIP-cons}. Straightforward CPLEX results}}
  \label{tb:3}
 \vskip .2cm
\begin{tabular}{l|rrrrrr}
  \hline
 &&&&&&\\[-1.5ex]
  $Inst$ &   $z_{LP}$&$t_{LP}$&$\underline{z}_{CPX}$&$\overline{z}_{CPX}$&$t_{\overline{z}_{CPX}}$&$GAP_{CPX}\%$\\[0.5ex]
\hline
 I1&7110.0&0.18&7467.3&7468.4&518&0.00\\
 I2&8749.2&0.48&9081.2&9240.6&43200 (39300)&1.72\\
 I3&8107.6&0.75&8283.6&8488.6&43200 (39260)&2.41\\
 I4&18367.2&1.95&18746.7&19712.2&43200 (38700)&5.04\\
 I5&12936.8&2.79&13034.1&13547.2&43200 (11975)&3.78\\
 I6&36031.6&14.54&36410.7&41452.2&43200 & 12.16\\
  I7&19849.5&49.70&19949.6&21539.1&43200&7.38\\
  I8&21609.9&96.10&21609.9&1.67224$e^8$&43200&99.99\\
\hline
\end{tabular}
\end{table}

\subsection{Matheuristic SCS4B results}
Table \ref{tb:4} provides the main results obtained when model LIP \eqref{LIP-obj}--\eqref{LIP-cons} is solved by the matheuristic SCS4B.
The headings are as follows:
$\underline{\underline{z}}$ and $t_{\underline{\underline{z}}}$, lower bound \eqref{lb-week} of the optimal solution value retrieved in Step 2 from the solution values of the strip and stack $c$-submodels \eqref{strip-c-model} and \eqref{stack-c-model}, respectively, and wall time;
$\underline{z}$ and $t_{\underline{z}}$, lower bound  \eqref{lb} of the optimal solution value retrieved in Step 2 from the solution values of the  $c$-submodels \eqref{c-model} and  wall time;
$\overline{z}$ and $tt_{\overline{z}}$, incumbent solution value retrieved in Steps 3 and 4 related to the first stage variable fixed values retrieved from the `best' solution of the $c$-submodels \eqref{c-model} and total
wall time, respectively, where the time $t_{\underline{z}}$ is also summed;
$solver=\{CPX, LSH\}$, where $CPX$ is straightforward CPLEX and $LSH$ is the local search heuristic (see \cite{cdap23})
for solving $\omega$-submodel \eqref{h-omega}, for
$\omega\in\Omega$, in Step 3 to obtain the incumbent solution, depending upon the dimensions of the submodel;
$\overline{\overline{z}}$ and $tt_{\overline{\overline{z}}}$,
incumbent solution value retrieved from the solution in Steps 3 and 4 related to the
`best' first stage variable fixed values from the strip and stack $c$-submodels \eqref{strip-c-model} and \eqref{stack-c-model}, respectively,
and total wall time;
$GAP\%$, optimality gap of the incumbent solution expressed as $100\cdot{{\overline{z}-\underline{z}'}\over {\overline{z}}}$,
where $\underline{z}'=\max\{\underline{z}, \, \underline{\underline{z}}\}$; and
$GR(\overline{z})$, goodness ratio of the SCS4B incumbent solution versus the one provided
by straightforward CPLEX, expressed as $\frac{\overline{z}}{\overline{z}_{CPX}}$.
Note that the parameter $\rho=0$ in all instances and thus $out=0$ (see Section \ref{sec:SCS4B});
this means that the trial $(\hat{\alpha},\hat{\beta})$ is not considered to obtain the incumbent solution if it implies at least one outsourcing scenario $\omega$.

\begin{figure}[h]
\begin{center}
\includegraphics[width=12.cm,viewport=93 529 504 741,clip=]{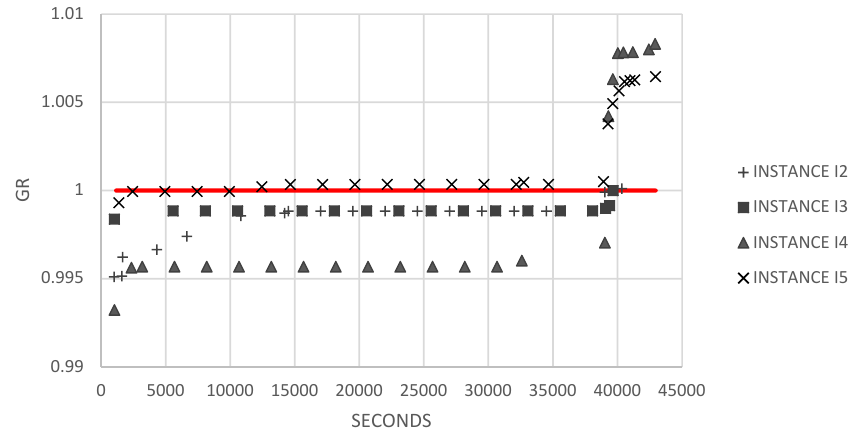}
\caption{SCS4B goodness ratio trajectory for instances I2 to I5}\label{fig:gr1}
\end{center}
\end{figure}

\begin{table}[!h]
\small
  \caption{\footnotesize{Model LIP \eqref{LIP-obj}--\eqref{LIP-cons}. SCS4B results}}
  \label{tb:4}
 \vskip .2cm
\begin{tabular}{l|rrrrcrrrrrr}
  \hline
 &&&&&&&&&&\\[-1.5ex]
 $Inst$ &  $\underline{\underline{z}}$&$t_{\underline{\underline{z}}}$& $\underline{z}$&$t_{\underline{z}}$&$solver$ &$\overline{z}$&$tt_{\overline{z}}$
 & $\overline{\overline{z}}$ & $tt_{\overline{\overline{z}}}$ & $GAP \%$ & $GR(\overline{z})$\\[0.5ex]
\hline
 I1 & 7138.0 & 2 & 7385.6 & 45 &$CPX$ & 7488.4 & 94 &  7509.0 & 5 &1.37 & 1.0028\\
 I2&8826.2&2&9057.8&163&$CPX$ &9241.6&363           &  9469.2& 13& 1.98&1.0001\\
 I3&8200.6&3&8184.0&426 & $CPX$&8488.6&557           &8594.5 &13 & 3.39&1.0000\\
 I4&18392.2&27&19296.8&1805& $CPX$ &19847.8&2070      &20080.2 &837 & 2.77&1.0069\\
 I5&12961.8&214&11808.0&1834&$CPX$ &13634.8&2820     & 13658.4 &568 & 4.93&1.0065\\
 I6&36243.6&3321&36110.3&18000&$LSH$ &36613.0&21069  & 37137.0& 5789& 1.01&0.8832\\
 I7&19728.2&9034&17864.5&21600&$CPX \& LSH$&20592.4&23990&20863.3&12236&4.19&0.9560\\
 I8&21610.0&25200&12292.8&25200&$CPX \& LSH$&23441.8&34108&23482.2&34525&7.80&$\approx 0$\\
\hline
\end{tabular}
\end{table}
It is worth pointing out that the lower bounds for the instances I3 and I5 to I8 are such that $\underline{\underline{z}}> \underline{z}$.
This is due to the fact that partitioning into a large
number of clusters was necessary for obtaining $\underline{z}$ \eqref{lb} in those instances, and additionally,
a compact model with all scenarios has been considered in the related strip and stack $c$-submodels \eqref{strip-c-model} and \eqref{stack-c-model}, respectively,
to obtain $\underline{\underline{z}}$ \eqref{lb-week}.
Moreover, it can be seen that the decomposition of the
problem by relaxing the relationships between the strip and stack doors to obtain stronger lower bounds \eqref{lb-week} for those instances
is not reflected in the strength of the corresponding feasible solutions (note that $\overline{\overline{z}}> \overline{z}$);
however, the wall time for obtaining $\underline{\underline{z}}$ is smaller than that for $\underline{z}$.
It is clearly shown that SCS4B outperformed straightforward CPLEX in terms of wall time,
it provides practically the same incumbent solution value for the instances I1 to I5 (see column GR), and
gives a smaller incumbent value for the large instances (I6 to I8).
Observe that SCSS4B outperforms straightforward CPLEX in terms of both the incumbent solution value and wall time for the large instances I6, I7, and I8;
as a matter of fact, straightforward CPLEX, from a practical point of view, does not provide any solution for the largest instance, I8 (see Table \ref{tb:3}).
On the other hand, both approaches provide practically the same solution value for instances I1 to I5.
Note the large amount of time that CPLEX requires to obtain a similar solution to the one given by SCS4B;
it is worth considering the goodness ratio's evolution for instances I2 to I5 as the wall time increases, as depicted in Fig. \ref{fig:gr1}.

\section{Conclusions and future research agenda}\label{sec:con}
A two-stage stochastic binary quadratic (BQ) formulation and related MILP have been introduced for the Cross-dock Door Design problem under uncertainty in the set of origin and destination nodes and, thus, in the number of product pallets that must traverse the cross-dock from origin nodes to destination nodes in a given time period and with the potential disruption of the door capacity.
Two schemes are proposed for obtaining strong lower bounds of the optimal solution value of this NP-hard combinatorial problem, as well as a matheuristic constructive algorithm called  SCS4B.
These schemes take advantage of the GAP structures that are embedded in the model; in particular, they are embedded in its seemingly most efficient Linearized Integer Programming reformulation (LIP).
It has been computationally proved that the overall SCS4B scheme requires a reasonable wall time that is much smaller than the wall time required by straightforward CPLEX, and it provides practically the same incumbent solution value for small and medium-sized instances and a better solution value for larger instances.

\subsection{Main items on the future research agenda for CDDP-TS}\label{sec:agenda}
\begin{itemize}
\item Embarrassing parallelism can be exploited, via a collaborative mode (see \cite{Aldasoro17}), to solve the scenario cluster submodels to obtain strong lower bounds, and it can be used to solve the scenario submodels to obtain efficient solutions.
\item A two-stage distributionally robust optimization BQ approach is to be considered for risk-neutral and risk-averse modes (second-order stochastic dominance functional).
\end{itemize}

\section*{Data Availability Statement}

The data that support the findings of this study are available from the corresponding author,  L.F. Escudero, upon reasonable request.

\section*{Disclosure Statement}

The authors report that there are no competing interests to declare.

\section*{Acknowledgments}
This research has been partially supported by the research projects
RTI2018-094269-B-I00 and PID2021-122640OB-I00  (L.F. Escudero) and
PID2019-104933GB-I00 from the Spanish Ministry of Science and
Innovation and Grupo de Investigaci\'on EOPT (IT-1494-22) from the Basque government (M.A. Gar\'in and A. Unzueta).

\end{document}